\documentclass[10pt]{article}
\usepackage{amsmath,amssymb,amsthm,graphics,graphicx}
\usepackage[T1]{fontenc}
\usepackage{multicol}
\usepackage{epstopdf}
\begin{document}
\vspace{5pt}

\def\lc{\left\lceil}   
\def\rc{\right\rceil}
\def\lf{\left\lfloor}   
\def\rf{\right\rfloor}

\begin{center}
\vspace{2pt}
    {\Large{Tight lower bounds for connected queen domination problems on the chessboard \footnote{International Conference on Advances in Applied Probability, Graph Theory and Fuzzy Mathematics
11-14 January, 2014,
}}
 }
\\
{Sneha S. Venkatesan \footnote{Computer Science Student, UCLA}\\}
{S. M. Venkatesan \footnote{Bangalore, India}\\
}
 \end{center}
 \vspace{4pt}

\begin{abstract}

.\\
1. We first show a lower bound of $2N/3 - 1$ for the connected minimum queen domination (or cover) problem on the $N * N$ chessboard - the upper bound is only 2 higher at most and is easy to show.\\
2. We then define the $k$-colored connected minimum queen domination, and extend the above proof to show a lower bound of $(2N -  k - 2)/3$, where the parameter $k$ can be increased to get decreasing lower bounds $LB (N, k)$ until one reaches the simple domination lower bound of  $\lf N / 2 \rf$. \\
3. We also discuss extensions of the connected domination problem and additional directions.\\

\textit{Keywords - chessboard problems, n-queens, queen domination number, connected domination, lower bounds, graph vertex cover, independent set, higher dimensional chessboards}

\end{abstract}

\section{Introduction}

The problem of placing the minimum number of queens on a $N * N$ hole-free chessboard so that they cover the entire board is a well-known problem with a history perhaps dating to the 19th century. Lower and upper bounds for this problem (best currently known values for these are around $n/2$ [1, 8]) were introduced in the $80’s$.  The \textit{connected queen domination} problem however asks that the placement of the queens in an optimum solution be such that each queen is covered by another queen (this can be viewed as asking that a spanning tree or a $k$-tree forest exist in a graph $G(Q)$ formed by adding an edge between every queen-pair that cover each other in the solution). It is  intuitive that the lower bound for this connected domination problem will be greater than for the normal domination problem since the coverage count of an empty/uncovered row, column, or diagonal by a queen will be potentially reduced by half when we insist that it be in the same row (column or diagonal) as another queen. We show here that this intuition is correct and that connected domination requires $2n/3$ - $k$ queens, where $k$ is the number of trees in the spanning forest. We also show a matching upper bound making this the first significant problem in this space to have tight bounds up to an additive constant. We are currently investigating extensions of the problem definition and the proof techniques to new natural and interesting variants of the connected domination problem on the 2D chessboard (these variant definitions but not the proofs also extend to general graphs, not necessarily chess-board based), and also to higher dimensions.

\section{Review of Previous Lower Bound for Normal Domination}

A.	\textit{The Problem and solution context} \\
We use and follow the descriptions from [1] (it was  assumed there that the queens can attack one another, but the lower bound there will be valid even if we insist that a queen not attack any other queen, or that every queen must be attacked by at least one other queen). We derive the lower bound in terms of a rectangular (assumed hole-free) $N * N$ chessboard.\\
\\
B.	\textit{The framework for the lower bound proof} \\
We review the 2-dimensional $N * N$ solution from [1] first to describe the simple idea, and second to set the context for our solutions which follow similar ideas.  Number the rows and columns from the bottom and the left respectively, and denote a square in column x and row y by $(x, y)$. A queen on $(x, y)$ covers all the squares on column x, on row y, and on the two diagonals passing through $(x, y)$. Given some placement of queens on such a board, we define an empty square as one without a queen placed on it, and an uncovered square as an empty square not covered by any queen. We also define an empty row (column or diagonal) to be one all of whose squares are empty, and an uncovered row (column or diagonal) to be one at least one of whose squares is uncovered. The main result from [1] is the following:\\

Result 1 [1]: At least  $ \lf N/2 \rf$ queens are needed to cover an $N * N$ board. \\

Proof sketch: Assume some \textit{minimum} queen-cover has already been specified on the board, and let  $(x_1+1)$ and $(N - x_2)$ be the smallest and the largest indexed empty columns respectively). Also, let $(y_1+1)$ and $(N-y_2)$ be the smallest and the largest-indexed empty rows respectively. These sentinel rows and columns determine the regions labeled $Q_1$ through $Q_8$ in Figure 1, any of which can be null, and the region $Q_9$. If $Q_9$ is less than $N/2 * N/2$ in size, then, by our choice, there are at least $N/2$ queens already on the board, proving the theorem. Hence assume $Q_9$ to be non-null. \\

We will without loss of generality let the \textit{symbol} $Q_i$ stand for the number of queens in the region $Q_i$, $1 <= i <= 9$. The implied annulus A of empty rows and columns surrounding region $Q_9$ containing $(4N - 4 - 2x_1- 2x_2- 2y_1- 2y_2)$ empty squares must then be covered by the $\sum_{i} Q_i $ queens on the board. Any queen in regions $Q_1$ through $Q_4$ can cover at most 2 squares of$A$  (through one diagonal). A queen in regions $Q_5$ through $Q_8$ can cover at most 6 squares of A (through 2 diagonals and one row or column), and a queen in region $Q_9$ can cover at most 8 squares of $A$ (through two diagonals and a row and a column). Therefore\\

\begin{figure}
\centering
 \includegraphics[width=150pt]{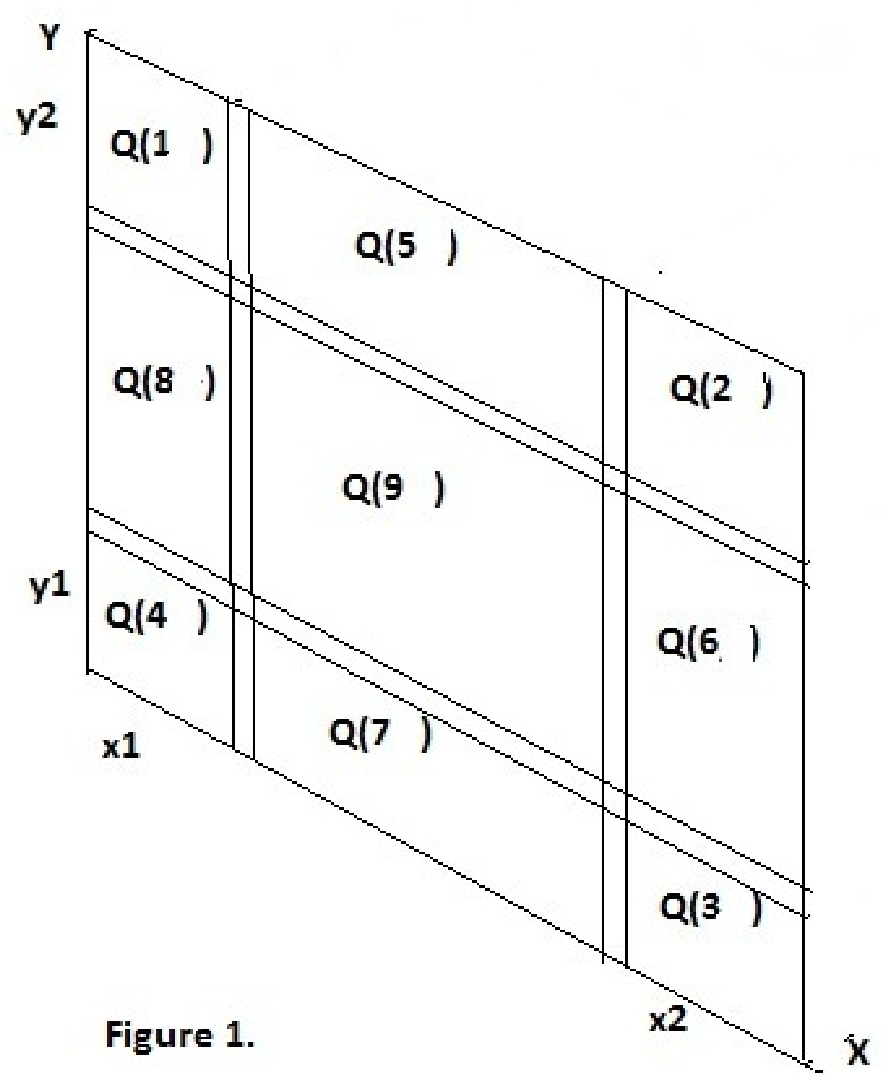} 
 \label{fig:1}
\end{figure}

$8.Q_9 + 6. \sum_{i=5}^{8} Q_i + 2. \sum_{i=1}^{4} Q_i >= 2(2N- 2-x_1-x_2-y_1-y_2)$\\
We add a slack of $2. \sum_{i=1}^{4} Q_i $ to the left hand side giving\\
$8.Q_9 + 6. \sum_{i=5}^8 Q_i + 4. \sum_{i=1}^{4} Q_i >= 2(2N- 2-x_1-x_2-y_1-y_2) $     (I) \\
\textit{This slack would give us room to deal with Subcase (b) discussed later.}
Counting outside the annulus by our construction, we get\\
$Q_3 + Q_4 + Q_7 >= y_1,  $\\
$Q_1 + Q_2 + Q_5 >= y_2,  $\\ 
$Q_4 + Q_1 + Q_8 >= x_1,  $\\
$Q_2 + Q_3 + Q_6 >= x_2   $\\
Summing and doubling gives\\
$2. \sum_{i=5}^{8} Q_i + 4. \sum_{i=1)}^4 Q_i >= 2(x_1+ x_2+y_1+ y_2) $       (II)\\
Combining these, we get  \\
 $\sum_{i=1}^{9} Q_i 	>= 1/8. \lc 4N - 4  \rc >= \lf N/2 \rf $      (III)

\section{Connected Domination}

Here we require that every queen on the board in a given solution cover at least one other queen. We shall rewrite the proof from the previous section as follows.\\

Let the quantity COMMON-D (COMMON-R or COMMON-C) denote the commonality of queens on a diagonal (row or column) - that is, if two queens share the diagonal (row or column), then the commonality added to COMMON-D (COMMON-R or COMMON-C) is 1, and if p queens share the diagonal (row or column), then the commonality added is p-1. Note that this is the same as the number of edges added to G(Q).\\

For better understanding, we will split COMMON-R as COMMON-R(I) and COMMON-R(II),  which stand for the effects on inequality (I) or inequality (II) respectively. We will split COMMON-C similarly into COMMON-C(I) and COMMON-C(II). We will also split COMMON-D as COMMON-D (on Annulus) and COMMON-D (not on Annulus). \\
  
First observe that our condition of \textit{two or more queens being in the same diagonal (row or column} : we use the term line to stand for any of these three terms) can be classified into these subcases: \\

a) When that diagonal crosses the annulus A: this will be accounted for by modifying inequality (I), and this will have effect on COMMON-D (on annulus).\\

b) When that diagonal does not cross the annulus (note they cross twice if they do): this is a special case which will be accounted for by modifying inequality (I), making use of the slack of $2. \sum_{i=1}^4 Q_i $ that was mentioned in the last section\\

c) When that row or column crosses the annulus A: this will be accounted for by modifying inequality (I), and this will have effect on COMMON-R(I) or COMMON-C(I).\\

d) When that row or column does not cross the annulus: will be accounted for by modifying inequality (II), and this will have effect on on COMMON-R(II) or COMMON-C(II).\\

In subcase (a), the influence of this diagonal on the annulus A in (I) is reduced by $2 * (p-1)$ where p is the number of queens in that diagonal.\\

We need to be very careful in our definition of when to add an edge to the G(Q), because p queens on a line will otherwise contribute a clique to G(Q) which will invalidate the determination of when G(Q) becomes connected based on edge cardinality - we shall add an edge between two queens  only when they \textit{directly} see each other or in other words no other queen is in between the two queens- therefore we add only $p-1$ edges for the line. With this definition of edge addition to G(Q), we would get a simple path of length $p-1$ connecting the queens that lie in the same row (column or diagonal). This ensures that the final spanning tree of G(Q) is embeddable in the same board (note that one of any two crossed edges that connect two disjoint paths causes a cycle and can be removed to construct a larger sub-tree from the two paths, thus growing the sub-tree until it becomes fully connected). The implications of this we will not go into in this paper.\\

In subcase (c), the influence of the row or column on the annulus A in (I) is again reduced by $2 * (p-1)$ with a similar definition of p.\\

In subcase (d), the row or column does not touch the annulus, that is, the effect is only on inequality (II). That is, for every such row or column shared by p queens, $2*(p-1)$ influence is charged to COMMON-R(II) or COMMON-C(II).

In subcase (b), there is no influence of the diagonal on the annulus A in (I), but this can happen only when the diagonal is in one of the four corners of the board without touching the annulus. For example, such a diagonal can be completely in $Q_1$  going from top-right to bottom-left, or cross from $Q_8$ to $Q_1$ to $Q_5$, and so on. Depending on the layout of the annulus (which we do not know before), the reduced effect on the annulus is still $2 * (p-1)$. Note that if the situation involves queens in $Q_5$ through $Q_8$ (but not $Q_1$) then that line must cross the annulus  and this does not cause any problems in the inequality. \textit{Therefore} we can uniformly treat the most important case of diagonal going through $Q_1$ on par with $Q_5$, since we have incremented the left hand side of (I) by the slack $2. \sum_{i=1}^4 Q_i$.\\

In subcase (b), the influence on the annulus of a queen in $Q_5$ through $Q_8$ can be 2, 4, or 6.  When it is only 4 (that is, one line does not cross the annulus) and the queen shares that line with another queen (if such a queen does not share the line with another queen, we get some beneficial slack in the inequality), then it adds 2 to the COMMON-D (not on Annulus) term. When it is only 2 (that is, two lines do not cross the annulus) and the queen shares either or both of these lines with another queen (if it does not share, we get beneficial slack as said before), then it adds 2 or 4 to COMMON-D (not on Annulus) term.\\
Therefore we will rewrite the inequalities as\\
$8.Q_9 + 6. \sum_{i=5}^8 Q_i + 4. \sum_{i=1}^4 Q_i $ \\
- 2*COMMON-R(I) - 2*COMMON-C(I) \\
- 2*COMMON-D (on annulus) -  2*COMMON-D (not on Annulus) \\
$>= 2(2N- 2-x_1-x_2-y_1-y_2) $          (IV)\\
$2. \sum_{i=5}^8 Q_i + 4. \sum_{i=1}^4 Q_i $
- 2*COMMON-R(II) - 2*COMMON-C(II) \\
$>= 2(x_1+ x_2+y_1+ y_2)$        (V) \\

Combining these, we get  \\
$ 8.\sum_{i=1}^9 Q_i  >= 4.N- 4 $ + \\
2* (COMMON-R(I) + COMMON-C(I) + COMMON-D (on annulus) +\\
 COMMON-D (not on Annulus) + COMMON-R(II) + COMMON-C(II)) \\
$8.\sum_{i=1}^9 Q_i  >= 4.N- 4 $ + \\
2* (COMMON-R + COMMON-C + COMMON-D)\\
The way we have defined, edges are added to $G(Q)$ only when a queen can directly see another queen (note that visibility in 2-D does not form a planar connection). Therefore \\
(COMMON-R + COMMON-C + COMMON-D) $>= (N_Q -1)$\\
(where $N_Q$ is the number of queens on the board) because that is the minimum number of edges needed in $G(Q)$ for connected domination. Since this gives $ (8.N_Q >= 4.N- 4 + 2* (N_Q -1))$, we get\\

RESULT 2: At least $\lc 2N/3 - 1 \rc $ queens are necessary for connected domination\\
Similar results can be shown for rectangular boards. \\
Upper Bounds: Near tightness of this lower bound can be proved by a simple construction. Consider a 3 X 3 almost equal block split of a N X N board (say  N is a multiple of 3, the other residuals can be handled easily as well as in [1]) such that the center-column of blocks has one less column, and the left-column of blocks has one more column. Then placing queens on the (bottom left to top right) diagonals in the bottom right and top left (in a right-aligned manner) blocks will cover all the squares except the first column: then we need one more queen in the leftmost column to cover this and one queen in the rightmost column to connect the two sets, giving an upper bound $ \lc 2N/3 + 1 \rc$.\\

\section{K-Colored Connected Domination}
Here we assume that there are k colors of queens available, and we require that all queens on the board of any given color form a connected set (in the sense described in the previous section). Then we can reformulate the result from last section as\\
(COMMON-R + COMMON-C + COMMON-D) $ >= (N_Q - k)$\\
Therefore, $ (8.N_Q >= 4.N- 4 + 2* (N_Q -k))$.  We then get\\

RESULT 3: At least $ \lc (2N - k - 2)/3 \rc $ queens are necessary for k-colored connected domination\\

If we parameterize the problem and the solution by increasing $k$, the lower bound of  $LB (N, k) = (2N - k - 2)/3 $ keeps falling, until one reaches the simple domination case, where
$LB (N) = ( 2.N - LB (N) - 2 )/3 $ 

Simplifying, we get  $ LB (N) >= (N - 1)/2$, which not surprisingly is the simple domination lower bound.\\

\section{Other Directions and Conclusions}
Previously we noted that we add edges to G(Q) only when a queen can directly see another queen. We can however add the other edges to G(Q) when more than 2 queens, say l, are on a single line, in which case l(l-1)/2 edges are added. We are investigating the very interesting question of what is the maximum number of queens a single row, column, or diagonal can have, \textit{without} deviating from the lower bound. \\

We are also investigating extending the connected domination arguments to a k-dimensional chessboard, but this is a much harder problem and is out of the scope of this work. It is also possible to define colored versions (as above) in many of the other chessboard problems, perhaps even in the notion of polynomials studied in some of them.\\

Because the geometry of chessboard domination breaks when one moves to connected domination in general graphs (which does have interesting applications), we believe these techniques cannot be employed there. \\

However, it looks like the double and triple counting arguments (as introduced in [1]) - leading to self-canceling sums of single or multiple types of objects (queens in the case of [1]) in the geometric combinatorial intersections of axis-different sub-arrays - will be of quite some use in proving lower bounds (and showing approximation algorithms and schemes) in computational geometry, orthogonal subspaces, algebra, lattices, weighted and unweighted packing and knapsack problems in 2- and 3-D, and other geometries as well. The ideas from [1] and from this paper can be recast in many different ways, but seem to be the only  way to visualize and investigate this and the new problems.\\

\textit{Note}: When this result was communicated to Dr. Weakley, he has informed us that a lower bound for the case k = 1 has been announced by P. Burchett in an abstract in the Thirty-Seventh Southeastern International Conference on Combinatorics, Graph Theory and Computing, 2006.\\

\end{document}